\documentclass{article}[12pt]
\usepackage{latexsym,amssymb,amsthm,amsfonts,amsmath,epsfig,psfrag
}
\newtheorem{lemma}{Lemma}
\newtheorem{theorem}{Theorem}
\newtheorem{prop}{Proposition}

\def \H{{\mathbb H}}

\begin{document}

\begin{center}
{\Large  \bf
A series of word-hyperbolic Coxeter groups
}

\vspace{15pt}
{\large Anna Felikson, Pavel Tumarkin}
\end{center}

\vspace{8pt}

\begin{center}
\parbox{11cm}%
{\small
{\it Abstract.}
For each positive integer $k$ we present an example of Coxeter system 
$(G_k,S_k)$ such that $G_k$ is a word-hyperbolic Coxeter group,  
for any two generating reflections $s,t\in S_k$ the product $st$ has finite 
order, and the Coxeter graph of $S_k$ has negative inertia index equal to $k$.
In particular, $G_k$ cannot be embedded into pseudo-orthogonal group  
$O(n,k-1)$ for any $n$ as a reflection group with generating set being 
reflections w.r.t. $(G_k,S_k)$.
}
\end{center}

\vspace{8pt}

\noindent
A {\it Coxeter system} $(G,S)$ consists of a Coxeter group $G$ and a 
generating set $S$ of $G$ with defining relations 
$s^2=t^2=(st)^{m_{st}}=1$ for some $m_{st}=2,3,\dots,\infty$ for any
pair $s,t\in S$. The set of 
{\it reflections} of the group $G$ with respect to $(G,S)$ is the set 
of all conjugates of elements of $S$. 
   
The purpose of this short note is to prove the following

\begin{theorem}
For each positive integer $k$ there exists a Coxeter system 
$(G_k,S_k)$ such that $G_k$ is a word-hyperbolic Coxeter group,
and $G_k$ cannot be embedded into pseudo-orthogonal group
$O(n,k-1)$ for any $n$ as a reflection group with generators being 
reflections with respect to $(G,S)$.

\end{theorem}

Below we give an explicit construction of these groups and Coxeter 
systems, write down
a presentation by standard generators and relations, and draw
Coxeter graphs of the groups.

The question was introduced by L.~Potyagailo while discussing the
paper~\cite{JS}.
In that paper, T.~Januszkiewicz and J.~Swiatkowski
constructed right-angled word-hyperbolic Coxeter groups of arbitrary
large virtual cohomological dimension. However, we have not managed to
understand possible signatures of these groups.

We are grateful to L.~Potyagailo for useful discussions and remarks. 
We also thank the Max-Planck-Institut f\"ur Mathematik in Bonn for 
hospitality.

\subsubsection*{1. Construction}
We refer to~\cite{B} for background on Coxeter groups, their
matrices and graphs.

Let $G_1$ be a group generated by 5 involutions $\{r_1,\dots,r_5\}=S_1$ and
relations\\
$\bullet\ (r_ir_j)^5=1$ \ if \ $(i,j)=(4,5)$;\\
$\bullet\ (r_ir_j)^3=1$ \ if \ $(i,j)=(1,2)$, $(2,3)$ or $(3,4)$;\\
$\bullet\ (r_ir_j)^2=1$ \ otherwise.

$(G_1,S_1)$ may be thought as a cocompact group of isometries of the
hyperbolic space $\H^4$. Its fundamental polytope is a simplex with
Coxeter diagram $\Sigma_1$:
\begin{center}
\psfrag{1}{\small $1$}
\psfrag{2}{\small $2$}
\psfrag{3}{\small $3$}
\psfrag{4}{\small $4$}
\psfrag{5}{\small $5$}
{\large $\Sigma_1$}\qquad\epsfig{file=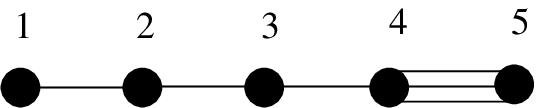,width=0.25\linewidth}
\end{center}

\begin{figure}
\begin{center}
\psfrag{1}{}
\psfrag{2}{}
\psfrag{3}{}
\psfrag{4}{}
\psfrag{5}{}
\epsfig{file=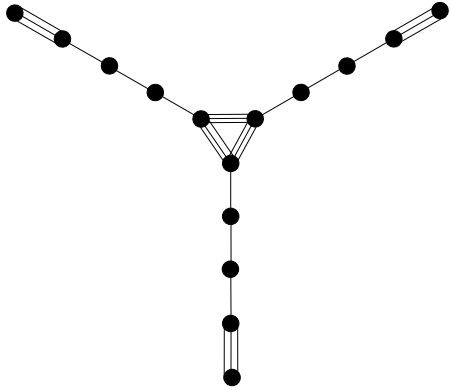,width=0.42\linewidth}
\end{center}
\caption{The Coxeter graph $\Sigma_3$ of the group $G_3$.}
\label{3}
\end{figure}

Now define a group $G_k$. It consists of $k$ copies of $G_1$ (an
$m$-th copy is generated by 5 involutions
$r_{5(m-1)+1},\dots,r_{5m}$)  with the following additional relations
between generators of different copies:\\
$\bullet\ (r_ir_j)^5=1$ \ if \ $i\equiv j\equiv 1\mod{5}$;\\
$\bullet\ (r_ir_j)^2=1$ \ otherwise.\\
Let $S_k$ be the set of generating reflections $\{r_1,\dots,r_{5k}\}$. 
The Coxeter graph $\Sigma_k$ of $(G_k,S_k)$ consists of $k$ copies of
Coxeter graph $\Sigma_1$ sticking out of a complete graph on $k$
vertices with triple edges. The graph $\Sigma_3$ is drawn on
Fig.~\ref{3}.

Since any two involutions $r_i$ and $r_j$  generate a finite group,
the signature of $\Sigma_k$ is uniquely defined.
To prove the theorem, we show that $G_k$ is word-hyperbolic and the
signature of $\Sigma_k$ equals $(4k,k)$. Then Lemma 12 of the 
paper~\cite{V71} implies that the set of reflections $S_k$ does not 
preserve any quadratic form with negative inertia index less than $k$. 
This means that $G_k$ cannot be embedded into $O(n,k-1)$ for any $n$ 
with all generators being conjugates of $\{r_i\}$. 

\subsubsection*{2. Main tools}
To prove that $G_k$ is word-hyperbolic, we use the following
result of G.~Moussong~\cite{M}.
\begin{prop}[\cite{M}, cor. of Th.~17.1]
\label{hyp}
Let $\Sigma$ be a Coxeter graph of a Coxeter group $G$ without
parabolic subgraphs. Then $G$ is word-hyperbolic if and only if no
pair of unjoined subgraphs $\Sigma'$, $\Sigma''$ of $\Sigma$
corresponds to two infinite groups.

\end{prop}

\noindent
To compute the determinant of $\Sigma_k$ we use the following result of
E.~Vinberg~\cite{V}.

The {\it weight} of an $(m-2)$-tuple edge of $\Sigma$ is a real number
equal to $\cos(\frac{\pi}{m})$. For any cyclic path $\gamma$ of  $\Sigma$
denote by $p(\gamma)$ the product of the weights of all edges of
$\gamma$.

\begin{prop}[\cite{V}, Prop. 11]
\label{det}
The determinant of a Coxeter graph $\Sigma$ is equal to the sum of
products of the form
$$(-1)^s p(\gamma_1)\dots p(\gamma_s),$$
where $\{\gamma_1,\dots,\gamma_s\}$ ranges over all unordered
collections (including the empty one) of pairwise disjoint oriented
cyclic paths of $\Sigma$.

\end{prop}

\subsubsection*{3. Proof of the theorem}
First, we prove that $G_k$ is word-hyperbolic. It is clear that any
subgraph of $\Sigma_k$ corresponding to infinite subgroup contains at
least one vertex of the complete subgraph on $k$ vertices. By
Prop.~\ref{hyp}, this implies word-hyperbolicity of $G_k$.

Hiding all vertices of the complete subgraph in $\Sigma_k$, we obtain
a graph corresponding to the finite group $kH_4$. So, the positive
inertia index of $\Sigma_k$ is at least $4k$. To prove that the
negative inertia index of $\Sigma_k$ equals $k$, it is sufficient to
show that the sign of ${\text{det}}(\Sigma_k)$ equals $(-1)^k$. We finish the
proof by the following

\begin{lemma}
$${\text{\em det}}(\Sigma_k)=\frac{(2-\sqrt{5})^k}{2^{5k}}(k+1)$$
\end{lemma}

\begin{proof}
Denote $\text{det}(\Sigma_k)$ by $d_k$, and denote by $d$ the determinant of
$H_4$. Prop.~\ref{det} implies the following recurrent formula for $d_k$:
$$d_{k+1}=d_kd_1+\sum\limits_{m=1}^{k}(-1)^m m!
\text{\normalsize ${k\choose m}$} (-\cos\,(\frac{\pi}{5}))^{m+1} d^{m+1} d_{k-m}.
\eqno{\text{($*$)}}$$
Here ${k\choose m}$ is the number of $(m+1)$-simplices with one fixed
vertex in the complete graph on $k$ vertices, and $m!$ is the number
of oriented Hamiltonian cycles in 1-skeleton of $(m+1)$-simplex.
Notice now that
$$d_1=\frac{2-\sqrt{5}}{16},\quad d=\frac{7-3\sqrt{5}}{32},\quad
-\cos\,(\frac{\pi}{5})=\frac{-1-\sqrt{5}}{4},\quad
-d\cos\,(\frac{\pi}{5})=\frac{2-\sqrt{5}}{32},$$
and assume that the lemma is true for all $l\le k$ (we put $d_0=1$,
so the lemma is true for $k=1$). Then the expression~{($*$)} implies that \\

\noindent
$d_{k+1}=\frac{(2-\sqrt{5})^k}{2^{5k}}(k+1)\frac{2-\sqrt{5}}{16}+
\sum\limits_{m=1}^{k}(-1)^m m!\text{\normalsize ${k\choose m}$}
\left(\frac{2-\sqrt{5}}{32}\right)^{m+1}
\frac{(2-\sqrt{5})^{k-m}}{2^{5(k-m)}}(k-m+1)=$\\
{}\hfill{
$=\frac{(2-\sqrt{5})^{k+1}}{2^{5(k+1)}}
\left(2(k+1)+\sum\limits_{m=1}^{k}(-1)^m
\frac{k!}{(k-m)!}(k-m+1)\right)=
\frac{(2-\sqrt{5})^{k+1}}{2^{5(k+1)}}(k+2).$}

\end{proof}

\vspace{35pt}
\noindent
Independent University of Moscow\\
e-mail:
\phantom{ow} felikson@mccme.ru\\
\smallskip
\phantom{owowwww} pasha@mccme.ru

\end{document}